\documentclass{ifacconf}

\usepackage{graphicx}      
\usepackage{amsmath,mathtools,amssymb}
\usepackage{amsfonts}
\usepackage{xcolor}
\usepackage{tikz}
\usepackage{pgfplots}
\usetikzlibrary{external}
\usepackage{epstopdf}
\usepackage[font=footnotesize, labelfont={sf,bf}, margin=1cm]{caption}
\usepackage{graphicx, import}
\usepackage{algorithm}
\usepackage{algpseudocode}
\usepackage{booktabs}
\usepackage{caption}
\usepackage{natbib}
\usepackage{comment}
\usepackage{enumitem}
\usepackage{accents}

\newif\ifcommentandrea
\commentandreatrue

\newcommand{\norm}[1]{\left\lVert#1\right\rVert}

\newcommand{\en}[0]{}
\begin{document}
\begin{frontmatter}

\title{Stability Analysis of Real-Time Methods for Equality Constrained NMPC} 

\thanks[footnoteinfo]{This research was supported by the German Federal 
Ministry for Economic Affairs and Energy (BMWi) via 
eco4wind (0324125B) and DyConPV (0324166B), by 
DFG via Research Unit FOR 2401 and by the EU via ITN-AWESCO
(642 682)}

\author[First]{Andrea Zanelli,}    
\author[Second]{Quoc Tran Dinh,}    
\author[First]{Moritz Diehl}

\address[First]{University of Freiburg, Germany (e-mail: \{andrea.zanelli, moritz.diehl\}@imtek.uni-freiburg.de)}
\address[Second]{University of North Carolina at Chapel Hill, North Carolina, USA (e-mail: quoctd@email.unc.edu)}


\begin{abstract}                         
In this paper, a proof of asymptotic stability for the combined system-optimizer dynamics 
associated with a class of real-time methods for equality constrained nonlinear model predictive control is presented. 
General Q-linearly convergent online optimization methods are considered 
 and asymptotic stability results are derived for the case where a single iteration of the optimizer is carried 
out per sampling time. 
In particular, it is shown that, if the underlying sampling time is sufficiently short, asymptotic stability 
can be guaranteed. The results constitute an extension to 
existing attractivity results for the well-known real-time iteration strategy.
\end{abstract}

\end{frontmatter}

\section{Introduction}
Nonlinear model predictive control (NMPC) is an optimization based control strategy
that relies on the solution of parametric nonlinear noncovex programs (NLP) in order 
to compute an implicit feedback policy. Due to the considerable computational burden associated 
with the solution of such NLPs, NMPC has first found application in fields where the sampling times 
are generally slow enough to carry out the required computations. In particular, since the 1970s, 
successful applications of NMPC have been reported in the process control industry \citep{Rawlings2017}.
\par 
In more recent years, due to the significant progress in the development of efficient algorithms and 
software implementations and due to the increasing computational power available on embedded control units, 
NMPC has gradually become a viable strategy for applications with much shorter sampling times. Among others, 
we report on recent applications such as \citep{Zanelli2019b}, \citep{Albin2017} and \citep{Besselmann2015}, where 
sampling times in the milli- and microsecond range are met.
\par
In order to alleviate the computational burden associated with NMPC, inexact approaches are often exploited that rely 
on the computation of approximate solutions to the underlying NLPs. 
The so-called real-time iteration (RTI) method proposed in \citep{Diehl2002} exploits a single iteration 
of a sequential quadratic programming (SQP) algorithm in order to compute an approximate solution of the 
current instance of the nonlinear program. By using this solution to warmstart the SQP algorithm at 
the next sampling time, it is possible to track an optimal solution and eventually converge to it, 
as the system's state is steered to a steady state. An attractivity proof for such an algorithm is 
derived in \citep{Diehl2007b} for a simplified setting where inequalities are not present or inactive 
in the entire region of attraction of the closed-loop system. Other real-time 
algorithms with stability guarantees are the relaxed-barrier anytime MPC for linear-quadratic problems \citep{Feller2017}, 
and the approach for general nonlinear systems in \citep{Graichen2010} that assumes a decrease over time of the cost function. 
Finally, in the recent paper by \cite{Liao-McPherson2019}, under rather general settings, stability is established 
with the requirement that a sufficiently large number of iterations are carried out per sampling time.

In the present paper, the results in \citep{Diehl2007b} are extended such that not only attractivity, but also stability of the combined 
system-optimizer dynamics can be guaranteed. 

\subsection{Notation}
Throughout the paper we will denote the Euclidean norm  by $\| \cdot \|$, when referring to vectors, 
and, with the same notation, to the spectral norm 
\vspace{0.1cm}
\begin{equation}
    \| A \| \vcentcolon = \sqrt{\lambda_{\max}\left(A^{\top}A\right)},
\end{equation}
when referring 
to a (real) matrix $A$. All 
vectors are column vectors and we denote the concatenation of two vectors 
by 
\vspace{0.1cm}
\begin{equation}
    (x,y)\vcentcolon=\begin{bmatrix}x \\ y\end{bmatrix}.
\end{equation}
\par
\vspace{0.1cm}
We denote the derivative (gradient) of any function by 
$\nabla f(x) = \frac{\partial f}{\partial x}(x)^{\top}$ and the Euclidean ball 
of radius $r$ centered at $x$ as 
\begin{equation}
    \mathcal{B}(x,r) \vcentcolon=\{ y \,\vcentcolon \,\| x - y \| \leq r\}.
\end{equation}
Finally, we denote the identity matrix by $\mathbb{I}$.
\vspace{0.4cm}
\section{Independent System and Optimizer Dynamics}
Consider the following continuous-time optimal control problem:
\begin{equation}\label{eq:nmpc}
\begin{aligned}
&\underset{\begin{subarray}{c}
    s(\cdot), u(\cdot)
\end{subarray}}{\min}	&&\int_{0}^{T_f} l(s(t), u(t)) \text{d}t + m(s(T_f))\\
                            &\quad \text{s.t.} 		   	&&s(0) - x = 0,\\
                            & 							&&\dot{s}(t) = \phi(s(t), u(t)), \quad  t \in [0,T_f], \\
\end{aligned}
\end{equation}
where $s \vcentcolon \mathbb{R} \rightarrow \mathbb{R}^{n_x}$ and $u \vcentcolon \mathbb{R} \rightarrow \mathbb{R}^{n_u}$ represent the 
state and input of a system, respectively, whose dynamics
are described by 
$\phi \vcentcolon \mathbb{R}^{n_x}\times\mathbb{R}^{n_u}\rightarrow \mathbb{R}^{n_x}$.
The functions $l\vcentcolon \mathbb{R}^{n_x}\times\mathbb{R}^{n_u}\rightarrow \mathbb{R}$
and $m\vcentcolon \mathbb{R}^{n_x}\rightarrow \mathbb{R}$ represent the Lagrange and Mayer 
cost terms, respectively. Finally, $x$ is a parameter describing the current state of the system and 
we assume, without loss of generality, that $\phi(0,0) =0$.

We will regard a discretized version of \eqref{eq:nmpc} obtained with \textit{some} 
discretization method (e.g. multiple shooting):
\begin{equation}\label{eq:compact_nlp}
    P(x)\vcentcolon 
    \quad \quad
    \begin{aligned}
    &\underset{y}{\min}	&&f(y) \\
    &\,\, \text{s.t.}  &&g(y) + Bx = 0,\\
    \end{aligned}
\end{equation}
where $y \in \mathbb{R}^n$ describes the primal variables of the discretized problem,
$f\vcentcolon \mathbb{R}^n \rightarrow \mathbb{R}$ 
and $g\vcentcolon \mathbb{R}^n\rightarrow \mathbb{R}^{n_g}$. The parameter $x$ enters the 
equality constraints through the linear map defined by the constant matrix $B \in \mathbb{R}^{n_g \times n_x}$.
\begin{assum}\label{assum:cont}
    The functions $f$ and $g$ are twice continuously differentiable and have bounded first and second order derivatives. 
\end{assum}
\par
Let $X \subseteq \mathbb{R}^{n_x}$ denote the set that contains all the 
possible values of $x$ such that $P(x)$ has at least one solution. We will 
restrict our attention to the set $X_{\bar{V}}\vcentcolon= \{ x \vcentcolon V(x) \leq \bar{V}\}$, with 
\begin{equation}
    V(x) \vcentcolon= f(\bar{y}(x)),
\end{equation}
and where $\bar{y}(x)$ solves $P(x)$. 
Let $\bar{u}(x)$ denote 
the feedback policy
\begin{equation}
    \bar{u}(x) \vcentcolon= M_{u,y} \bar{y}(x)
\end{equation}
implicitly defined by $P(x)$, for some constant projection matrix $M_{u,y}$, where $\| M_{u,y} \| = 1$ is assumed 
for simplicity.
\subsection{System Dynamics} 
The system under control obeys the following sampled-feedback 
closed-loop dynamics:
\begin{defn}[System Dynamics]\label{defn:ts}
    Let the following differential equation describe the dynamics of the system 
    controlled using a constant input $u_0$: 
    \begin{equation}
        \begin{aligned}
            &\frac{d\psi}{dt}(t;x_0, u_0) = \phi(\psi(t; x_0, u_0), u_0), \\
            &\psi(0; x_0, u_0) = x_0. 
    \end{aligned}
    \end{equation}
    Here $\psi \vcentcolon \mathbb{R} \times \mathbb{R}^{n_x} \times \mathbb{R}^{n_u} \rightarrow \mathbb{R}^{n_x}$ describes the 
    trajectories of the system, $x_0$ denotes the state of the system at a given sampling instant and 
    $u_0$ the corresponding constant input.
    We will refer to the strictly positive parameter $T > 0$ as the \textit{sampling time} associated with the corresponding discrete-time system 
    \begin{equation}
        x_{\mathrm{next}} = \psi(T; x, u).
    \end{equation}
\end{defn}

In the following, we summarize an adapted version of standard assumptions used to ensure the stability
properties of the nominal NMPC scheme. 

\begin{assum}[Lyapunov Stability]\label{assum:lyapunov_stability}
    Assume that there exists positive constants $a_1,\, a_2, a_3$ and $T_0$ such that the following holds 
    for any $x \in X_{\bar{V}}$ and any $T \leq T_0$:
    \begin{subequations}
        \begin{align}
            a_1 \|x\|^2 \leq V(x) &\leq a_2 \|x\|^2, \label{eq:assum_lyap_1} \\ 
            V(\psi(T; x, \bar{u}(x))) - V(x) &\leq -T \cdot a_3 \|x\|^2. \label{eq:assum_lyap_2}
        \end{align}
    \end{subequations}
\end{assum}
\begin{rem}\label{rem:d2c_lyap}
    Notice that Assumption \ref{assum:lyapunov_stability}, for a fixed $T$ boils down to 
    the standard assumption for exponential asymptotic stability (see e.g. Theorem 2.21 in \citep{Rawlings2017}). Moreover, 
    the dependency on $T$ in \eqref{eq:assum_lyap_2} can be justified, for example, by assuming that a 
    continuous-time Lyapunov function $V_c(x(t))$ exists such that $\frac{\text{d}}{\text{d}t} V_c(x(t)) \leq -\underline{a} \| x \|^2$, for some 
    positive constant $\underline{a}$ and that $V(x)$ is a sufficiently good approximation of $V_c(x)$ in the following sense. 

    In particular, regard the simpler case in which the system under consideration is linear time-invariant, 
    i.e. $\dot{x}(t) = \phi(x(t),u(t)) = A_c x(t) + B_c u(t)$. Its discretized counterpart reads
    $x_{\text{next}} = A_d x + B_d u$, where
    \begin{equation*}
        A_d \vcentcolon= \exp{(A_c T_d)}, \, B_d \vcentcolon= \left(\int_{0}^{T_d}\exp{(A_c\tau) \text{d}\tau}\right)B_c.
    \end{equation*}
    When controlling a discrete-time linear time-invariant system with the linear feedback policy $u = K_dx$, 
    we know that, if $x^{\top}P_dx$ is a Lyapunov function for the resulting closed-loop system $x_{\text{next}} = (A_d + B_d K_d)x$, then
    it must satisfy the following discrete-time Lyapunov equation:
    \begin{equation}
        (A_d + B_d K_d)^{\top}P_d(A_d + B_d K_d) - P_d + Q_d = 0,
    \end{equation}
    for some positive-definite $Q_d$. It is easy to show that, if the discretization time $T_d$ is sufficiently small, then
    $x^{\top}P_dx$ is a Lyapunov function for the  continuous-time closed-loop system  $\dot{x}(t) = (A_c + B_c K_d)x(t)$, where we use 
    the discrete-time gain $K_d$. In particular, it suffices to show that a positive-definite matrix $Q_c$ exists such that the 
    following continuous-time Lyapunov equation is satisfied:
    \begin{equation}\label{eq:clyap_rem}
        (A_c + B_c K_d)^{\top} P_d + P_d (A_c + B_c K_d) + Q_c = 0.
    \end{equation}
    To this end, we note that $A_d = \mathbb{I} + T_d A_c + O(T_d^2)$ and $B_d = T_d B_c + O(T_d^2)$, such that we obtain 
    \begin{equation*}
        \begin{aligned}
            &&\bigg( \frac{A_d - \mathbb{I}}{T_d}& + O(T_d) + \left( \frac{B_d}{T_d} + O(T_d)\right)K_d\bigg)^{\top} P_d \\
            && &+ P_d\left( \frac{A_d - \mathbb{I}}{T_d} + O(T_d) + \left( \frac{B_d}{T_d} + O(T_d)\right)K_d\right) \\
            && &+ Q_c = 0
        \end{aligned}
    \end{equation*}
    and, multiplying by $T_d$,
    \begin{equation}
        \begin{aligned}
            &&( A_d - \mathbb{I} &+ B_d K_d)^{\top} P_d + P_d( A_d - \mathbb{I} + B_d K_d) \\
            && &+ E^{\top}E= -T_d Q_c,
        \end{aligned}
    \end{equation}
    where $E = O(T_d)$. Let $\tilde{A}_d \vcentcolon= A_d + B_d K_d$. Simplifying, 
    we obtain 
    \begin{equation*}
        \begin{aligned}
            &&( A_d - \mathbb{I} &+ B_d K_d)^{\top} P_d + P_d( A_d - \mathbb{I} + B_d K_d) + E^{\top}E\\
            && =&\tilde{A}_d^{\top} P_d + P_d \tilde{A}_d - 2 P_d + E^{\top}E\\
            && \preceq& \tilde{A}_d^{\top} P_d + P_d\tilde{A}_d  -2P_d + E^{\top}E\, + \\
            && &(\tilde{A}_d - \mathbb{I})^{\top}P_d(\tilde{A}_d - \mathbb{I})  \\
            && =&\tilde{A}_d^{\top}P_d\tilde{A}_d - P_d + E^{\top}E = - Q_d + E^{\top}E,
        \end{aligned}
    \end{equation*}
    where we have exploited the fact that $(\tilde{A}_d - \mathbb{I})^{\top}P_d(\tilde{A}_d - \mathbb{I}) \succeq 0$. Due to the fact 
    that $E^{\top}E = O(T_d^2)$ and $Q_d \succ 0$, we obtain that $-T_d Q_c \preceq -Q_d + O(T_d^2)$, and, for any sufficiently small discretization time $T_d$, there must exist a positive-definite $Q_c$ such that the continuous-time Lyapunov equation \eqref{eq:clyap_rem} is satisfied.
    Finally, with similar arguments it is possible to show that, if $T_d$ is small enough, for any sufficiently small 
    sampling time $T$, $x^{\top}P_dx$ is a valid Lyapunov function for the closed-loop system $x_{\text{next}} = (A_s + B_s K_d)x$, 
    where
    \begin{equation*}
        A_s \vcentcolon= \exp{(A_c T)}, \, B_s \vcentcolon= \left(\int_{0}^{T}\exp{(A_c\tau) \text{d}\tau}\right)B_c.
    \end{equation*}
    \color{black}

\end{rem}
\begin{assum}[Second Order Growth]\label{assum:sog}
    Assume that, for any $u \in \mathbb{R}^{n_u}$,
    for any $x \in X_{\bar{V}}$ and any $T \leq T_0$, the following holds: 
    \begin{equation}
        \begin{aligned}
            &&V(\psi(T;\, &x, u)) - V(x) \leq \\
            && &-T \cdot a_3 \|x\|^2 + T \cdot O(\| u - \bar{u}(x)\|^2). \label{eq:assum_lyap_2}
        \end{aligned}
    \end{equation}
\end{assum}
\begin{rem}\label{rem:sog}
Assumption \ref{assum:sog} can be informally justified by analyzing the properties of an underlying 
continuous-time Lyapunov function $V_c(x)$ and using an argument similar to the one used in Assumption 2.18 in \citep{Diehl2007b} 
in a discrete-time setting. In particular, under suitable differentiability assumptions, for any $\delta > 0$, we can write
\begin{equation*}
    \begin{aligned}
        &&V_c(x) &= \min_{u_{\delta}(\cdot)} \bigg\{\int_0^{\delta} l(\psi(\tau,x,u_{\delta}(\tau)), u_{\delta}(\tau))\text{d}\tau \\
        && &+ \tilde{V}_c(\psi(\delta, x, u_{\delta}(\cdot))) \bigg\} \\
        && &=\tilde{V}_c(\psi(\delta, x, u_{\delta}(\cdot))) + \!\! \int_0^{\delta} l(\psi(\tau,x,u_{\delta}(\tau)), u_{\delta}(\tau))\text{d}\tau \\
        && &+ \delta \cdot O (\| u_{\delta}(\cdot) - \bar{u}_{\delta}(\cdot)\|^2),
    \end{aligned}
\end{equation*}
where $\tilde{V}(x)$ is the optimal value function for a problem with shrunk horizon $T_f - \delta$.
Using the fact that $V_c(x) \leq \tilde{V}_c(x)$ and a quadratic lower bound on $l$ we can conclude
\begin{equation*}
    \begin{aligned}
        &&V_c(\psi(\delta, x, u_{\delta}(\cdot))) &\leq \tilde{V}_c(\psi(\delta, x, u_{\delta}(\cdot))) \\
        && \leq V_c(x) - \int_0^{\delta} &l(\psi(\tau,x,u_{\delta}(\tau)), u_{\delta}(\tau))\text{d}\tau \\
        && + \,\delta \cdot O (\| u_{\delta}&(\cdot) - \bar{u}_{\delta}(\cdot)\|^2) \\
        && \leq V_c(x) - \delta \cdot &\tilde{a}_3 \| x \|^2  + \delta \cdot O (\| u_{\delta}(\cdot) - \bar{u}_{\delta}(\cdot)\|^2),
    \end{aligned}
\end{equation*}
which justifies Assumption \ref{assum:sog}, if $V(x)$ is a ``sufficiently'' good approximation of $V_c(x)$ in the sense of Remark \ref{rem:d2c_lyap}.
\par
\end{rem}
The following slightly tailored version of \citep[Theorem 2.21]{Rawlings2017} provides
asymptotic stability of the closed-loop dynamics.
\begin{thm}\label{thm:conv_stability}
    Let Assumption \ref{assum:lyapunov_stability} hold. Then,
    the origin is an exponentially asymptotically stable equilibrium for the closed-loop
    system $x_{\text{next}} = \psi(T;x, M_{u,y}\bar{y}(x))$ for any $T \leq T_0$.
    \begin{pf}
        See \citep{Rawlings2017}.
    \end{pf}
\end{thm}
\vspace{-0.11cm}
\subsection{Optimizer Dynamics} 
The first-order necessary optimality conditions associated 
with \eqref{eq:compact_nlp} read as follows:
\begin{equation}\label{eq:kkt_int}
\begin{aligned}
    &0 && = \nabla f(y) + \nabla g(y)\lambda, \\
    &0 && = g(y) + Bx
\end{aligned}
\end{equation}
where $\lambda \in \mathbb{R}^{n_g}$ is the Lagrange multiplier associated with 
the equality constraints. Introducing  
\begin{equation}
    F(z) \vcentcolon = 
    \begin{pmatrix}
    \nabla f(y) + \nabla g(y)\lambda \\
    g(y)  
    \end{pmatrix},
\end{equation}
where $z = (y, \lambda)$, Equations \eqref{eq:kkt_int} can be expressed as 
\begin{equation}\label{eq:gen_kkt}
    0 = F(z) + Cx,
\end{equation}
where $C \vcentcolon= [\,\,\,0\,\, B^{\top}]^{\top}$.
Let $\bar{Z}(x)$ be the set of all stationary points satisfying \eqref{eq:gen_kkt} for a given $x$. 
The following assumptions are made.
\begin{assum}[Regularity]\label{assum:sreg}
    Assume that, for any $x \in X_{\bar{V}}$ there exists a unique solution $\bar{z}(x)$ and that
    second order sufficient conditions hold at $\bar{z}(x)$. 
    Moreover, assume that the steady-state solution is $\bar{z}(0) = 0$, i.e. $\bar{Z}(0) = \{0\}$.
\end{assum}
\begin{rem}\label{rem:sreg}
Although most classical NMPC stability results rely on the fact that the optimizer finds the global solution,
the assumption that a unique solution exists in $X_{\bar{V}}$ is somewhat restrictive and 
deserves further discussion.
On the one hand, this is similar to \citep[Assumption 2.3]{Diehl2007b} and to some extent a strong assumption. 
On the other hand, notice that, loosely speaking, it would suffice to assume that there is a unique ``branch'' of the 
solution manifold in a neighborhood $\mathcal{N}$ of the origin in $\mathbb{R}^{n_z} \times \mathbb{R}^{n_z}$. For 
the sake of simplicity, in this paper, we will however assume that $\bar{z}(x)$ is the unique solution in $X_{\bar{V}}$.
One further additional implication is the fact that the value function $V(x)$ is continuous.
This is again similar to what assumed in \citep{Diehl2007b} and in some of the work on inherent robustness of 
NMPC (see e.g. \citep{Pannocchia2011}). 
\end{rem}
\begin{prop}
    Let Assumptions \ref{assum:cont} 
    and \ref{assum:sreg} hold. Then there exist strictly positive constants
    $\sigma$, $\bar{r}_z$ and $\bar{r}_x$, such that, for any $x \in X_{\bar{V}}$, the solution $\bar{z}(x)$ is uniquely defined 
    over $\mathcal{B}(\bar{z}(x), \bar{r}_z)$ and the following holds:
    \begin{equation}
    \norm{\bar{z}(x'') - \bar{z}(x')} \leq \sigma \norm{x'' - x'},
    \end{equation}
    for any $x'', x' \in \mathcal{B}(x, \bar{r}_x)$.
    \begin{pf}
        The result is a direct consequence of the implicit function theorem 
         (also known as Dini's theorem).
    \end{pf}
\end{prop}
\begin{defn}[Optimizer Dynamics]\label{defn:od}
    Let the following \\discrete-time system describe the dynamics of the optimizer
    used to solve the parametric problem \eqref{eq:compact_nlp}
    \begin{equation}\label{eq:od}
        z_+ = \varphi(\psi(T; x, M_{u,z}z), z), \\
    \end{equation}
    where $\varphi \vcentcolon \mathbb{R}^{n_x} \times \mathbb{R}^{n_z} \rightarrow \mathbb{R}^{n_x}$ 
    and where $M_{u,z}$ is a properly defined projection matrix. For simplicity, 
    we will assume that $\norm{M_{u,z}} = 1$.
\end{defn}
\begin{assum}[Contraction]\label{assum:contraction_quad}
    There exists a radius $\hat{r}_z > 0$ and a positive constant $\hat{\kappa} < 1$ such that, for any given stationary point $\bar{z}(x)$ at $x\in X_{\bar{V}}$, 
    and any $z$ in $\mathcal{B}(\bar{z}(x), \hat{r}_{z})$, the optimization routine produces $z_+$ such that
    \begin{equation}\label{eq:contraction}
        \norm{z_+ - \bar{z}(x)} \leq \hat{\kappa} \norm{z - \bar{z}(x)}.
    \end{equation}
\end{assum}
Since we are interested in a real-time strategy that seeks 
an approximate solution to \eqref{eq:gen_kkt} as the parameter 
$x$ changes over time, we will exploit the following result on 
general real-time methods.
\begin{lem}\label{lem:contraction_quad}
    Let Assumptions \ref{assum:cont}, 
    \ref{assum:sreg} and \ref{assum:contraction_quad} hold. Then there exist strictly positive 
    constants $r_z$ and $r_x$, and finite positive constants $\sigma, \hat{\kappa} > 0$, with $\hat{\kappa} < 1$, such that, 
    for any $x$ in $X_{\bar{V}}$, any $z$ in $\mathcal{B}(\bar{z}, r_z)$, and 
    any $x_+$ in $\mathcal{B}_{\en}(x, r_x)$, it holds that
    \begin{align}\label{eq:track_contraction_quad}
     {\!\!\!\!}\begin{split}
         \norm{z_+ - \bar{z}(x_+)} & \leq  \hat{\kappa} \norm{z - \bar{z}(x)} +  \sigma \hat{\kappa} \norm{x_+ \!\!-\! x}_{\en}.
        \end{split} {\!\!\!\!}
    \end{align}
    \begin{pf}
        See e.g. \citep{Zanelli2019a}.
    \end{pf}
\end{lem}
\section{System-Optimizer Dynamics}\label{sec:gen_rti}
Theorem \ref{thm:conv_stability} and Lemma~\ref{lem:contraction_quad} provide 
key properties of the system and optimizer dynamics, respectively. In this section, 
we analyze the interaction between these two dynamical systems. 
In particular, we interested in analyzing the behavior of the following coupled system-optimizer dynamics.
\begin{defn}[System-Optimizer Dynamics]\label{defn:sod}
    Let the following equations describe the evolution of the state of the system and 
    of the optimizer's iterates:
    \begin{equation}\label{eq:so_dyn}
        \begin{aligned}
            &&x_+ &= \psi(T; x, M_{u,z}z), \\
            &&z_+ &= \varphi(\psi(T; x, M_{u,z}z), z). \\
        \end{aligned}
    \end{equation}
    In particular, after a simple coordinate change, we can refer to the following \textit{error} dynamics in compact form:
    \begin{equation}\label{eq:so_dyn_compact}
        \xi_+ = \Phi(T;\xi),
    \end{equation}
    where $\xi \vcentcolon = (x,z-\bar{z}(x))$ and $\Phi \vcentcolon \mathbb{R} \times \mathbb{R}^{n_x+n_z} \rightarrow \mathbb{R}^{n_x + n_z}$.
    We will refer to \eqref{eq:so_dyn_compact} as \textit{system-optimizer dynamics}.
\end{defn}
An attractivity proof
for real-time iterations based on the assumption that no inequalities are present in the 
problem formulation or that, equivalently, no active constraints are present in the 
region of interest, has been proposed in \citep{Diehl2007b}. Similarly, an attractivity 
proof where shifted iterations and a zero terminal constraint are used is derived  
in \citep{Diehl2005b}. 
\par

In  this paper, we prove instead asymptotic stability for the system-optimizer dynamics, which is in
general not implied by attractivity.
Moreover, the general contraction considered in Assumption 
\ref{assum:contraction_quad} covers a broader class of 
algorithms where the iterations need not be iterations of a sequential
quadratic programming method.

\subsection{Perturbed Error Contraction}

In order to be able to use the contraction from Lemma \ref{lem:contraction_quad}, we will make 
a general assumption on the behavior of the closed-loop system in a neighborhood of the equilibrium 
and for a bounded value of the numerical error.
\begin{assum}\label{assum:state_evolution_bound}
    There exist positive finite constants $L_{\psi,x}$ and $L_{\psi,u}$ such that,  
    for all $x \in X_{\bar{V}}$ and all $z$ such that $\| z - \bar{z}(x)\| \leq r_z$, the following inequality holds:
    \begin{equation}
        \| \psi(x, u) - x \| \leq T (L_{\psi,x}\| x \| + L_{\psi,u} \| M_{u,z}z \|).   
    \end{equation}
\end{assum}
\begin{prop}\label{prop:state_evolution_bound}
    Let Assumptions \ref{assum:cont}, 
    \ref{assum:sreg} and \ref{assum:state_evolution_bound} hold and define the following constants: 
    \begin{equation}
        \eta \vcentcolon = L_{\psi,u} + L_{\psi,x}\sigma 
    \end{equation}
    and
    \begin{equation}
        \theta \vcentcolon = L_{\psi,u}.
    \end{equation}
    Then, the following holds:
	\begin{equation}\label{eq:state_evolution_bound}
        \|\psi(x, u) - x \| \leq T (\eta\| x \| + \theta \| z - \bar{z}(x)\|),
	\end{equation}
	for any $x \in X_{\bar{V}}$ and any $z$ such that $\norm{z - \bar{z}(x)} \leq r_z$. 
    \begin{pf}
    In the following let $u = M_{u,z}z$ for ease of notation.
    Due to Assumption~\ref{assum:state_evolution_bound} we have that 
    \begin{equation}
        \| \psi(x, u) - x \| \leq T (L_{\psi,x}\| x \| + L_{\psi,u} \| u \|)   
    \end{equation}
    and, due to the assumption of regularity at the solution $\bar{z}(x)$ 
    and the fact that $\bar{z}(0) = 0$, we can write 
    \begin{equation*}
        \norm{u} \leq \norm{\bar{z}(x)} + \norm{z - \bar{z}} \leq \sigma \norm{x} + \norm{z - \bar{z}},
    \end{equation*}
    and the following holds:
    \begin{equation*}
        \begin{aligned}
            &&\norm{\psi(x,u) - x} \leq &\,\,T(L_{\psi,x} + L_{\psi,u}\sigma) \norm{x} \\
            &&                          &+ TL_{\psi,u}\norm{z - \bar{z}}. \qed
        \end{aligned}
    \end{equation*}
\end{pf}
\end{prop}
\begin{prop}\label{prop:contraction_x}
    Let Assumptions \ref{assum:cont}, 
    \ref{assum:sreg}, \ref{assum:contraction_quad} and \ref{assum:state_evolution_bound} 
hold. Moreover, define 
    \begin{equation}
        T'_1 \vcentcolon= \min \bigg\{\frac{r_x}{\eta r_{\bar{V}} + \theta r_z}, \frac{r_z(1 - \hat{\kappa})}{\sigma \hat{\kappa}(\theta r_z + \eta r_{\bar{V}})}\bigg\},        
    \end{equation}
    where
    \begin{equation}
        r_{\bar{V}} \vcentcolon = \sqrt{\frac{\bar{V}}{a_1}}.
    \end{equation}
    Then, for any $x \in X_{\bar{V}}$, any $\| z - \bar{z}(x)\| \leq r_z$ and 
    any $T \leq T_1 \vcentcolon= \big\{ T'_1, T_0\big \}$, the following holds:
\begin{equation}\label{eq:contr_z}
\| z_+ - \bar{z}(x_+)\|  \leq \, \kappa\| z - \bar{z}\| + T \cdot \gamma \| x \|, 
\end{equation}
where
\begin{equation}
    \kappa \vcentcolon = \hat{\kappa}(1+T \sigma\theta) < 1, \quad \gamma \vcentcolon = \sigma\hat{\kappa}\eta.
\end{equation} 
Moreover, $\| z_+ - \bar{z}(x_+)\| \leq r_z$.
\begin{pf}
    Given that $\| z - \bar{z} \| \leq r_z$ and that, due to Assumption \ref{assum:state_evolution_bound} and 
the definition of $T'_1$, we have $\| x_+ - x \| \leq r_x$ for all $x \in X_{\bar{V}}$ and we can apply the contraction from 
    Lemma \ref{lem:contraction_quad}: 
    \begin{align}
     {\!\!\!\!}\begin{split}
         \norm{z_+ - \bar{z}(x_+)}_{\en} & \leq  \hat{\kappa} \norm{z - \bar{z}}_{\en} +  \sigma \hat{\kappa} \norm{x_+- x}_{\en}.
        \end{split} {\!\!\!\!}
    \end{align}
    Applying the inequality from Proposition \ref{prop:state_evolution_bound}, we obtain
    \begin{equation}
    \| z_+ - \bar{z}(x_+)\|  \leq \, \kappa\| z - \bar{z}\| + T \gamma \| x \|, 
    \end{equation}
    where
    \begin{equation}
        \kappa \vcentcolon = \hat{\kappa}(1+T\sigma\theta), \quad \gamma \vcentcolon = \sigma\hat{\kappa}\eta.
    \end{equation} 
    \par
    Finally, due to the second term in the definition of $T'_1$, we have that $\| z_+ - \bar{z}(x_+) \| \leq r_z$ and $\kappa < 1$ since
    \begin{equation}
        \begin{aligned}
            &&T &\leq T'_1 = \frac{r_z(1 - \hat{\kappa})}{\sigma \hat{\kappa}(\theta r_z + \eta r_{\bar{V}})}  \\
            && & = \frac{r_z}{r_z} \frac{(1 - \hat{\kappa})}{\sigma \hat{\kappa}(\theta + \eta r_{\bar{V}}/r_z)} < \frac{1 -\hat{\kappa}}{\hat{\kappa}\sigma\theta}. \qed 
    \end{aligned}
    \end{equation}
\end{pf}
\end{prop}
\subsection{Perturbed Lyapunov Contraction}

In the following, we analyze the impact of the fact that the approximate 
feedback policy $M_{u,z}z$ is used, instead of the optimal one $M_{u,z}\bar{z}(z)$,
on the nominal Lyapunov contraction.
Throughout the rest of the paper, we will make use of the following shorthands:
\begin{equation}
V \vcentcolon= V(x)
\end{equation}
and
\begin{equation}
    V_+ \vcentcolon= V(\psi(T; x, M_{u,z}z))
\end{equation}
to denote the values taken by the optimal cost at the ``current'' state and 
at the one reached applying the suboptimal control action $M_{u,z}z$ starting 
from $x$. Similarly, we introduce 
\begin{equation}
E\vcentcolon = \|z - \bar{z}(x)\|
\end{equation}
and 
\begin{equation}
    E_+ \vcentcolon=\| \varphi(\psi(T; x, M_{u,z}z), z) - \bar{z}(\psi(T; x, M_{u,z}z))\|
\end{equation}
to denote the numerical error attained at the ``current'' and next iteration of 
the optimizer, where the error is computed with respect to the exact solution associated 
with the ``current'' and next state of the system. Note that, although $V, V_+, E, E_+$ depend on 
$(x,z)$, we omit that dependency for a more compact notation.

In Lemma \ref{lem:invariance_q2}, we will show that, under the condition that the sampling time $T$ is sufficiently small, 
we can guarantee positive invariance of a properly defined set.
\par
\begin{prop}\label{prop:inexact_lyapunov_contraction_q2}
    Let Assumptions \ref{assum:cont}, 
    \ref{assum:lyapunov_stability}, \ref{assum:sog} and \ref{assum:sreg} hold. Then,
    there exist finite positive constants $\mu, \bar{V}_q \leq \bar{V}$ and $r_q \leq r_z$, such that, for any $E \leq r_q$
    and any $x$ in $X_{\bar{V}_q}$, where $X_{\bar{V}_q} \vcentcolon = \{ x \vcentcolon V(x) \leq \bar{V}_q\}$, the following holds: 
    \begin{equation}\label{eq:coupled_contraction}
        V_+  \leq (1 - T \bar{a}) V + T\mu E^2, 
    \end{equation}
    with $\bar{a} \vcentcolon = \frac{a_3}{a_2}$.
    \begin{pf}
        Assumption \ref{assum:sog} implies that there must exist strictly positive constants $\mu, \bar{V}_q \leq \bar{V}$ and $r_q \leq r_z$
        such that the following holds
    \begin{equation}
        \begin{aligned}
            &&V(\psi(x, M_{u,z}z)) &\leq V(x) - T a_3 \| x \|^{2} - T \mu E^2 \\
            &&           &\leq V(x) - T \frac{a_3}{a_2} V(x)  - T \mu E^2 \\ 
            &&           &=(1 - T \bar{a}) V(x) - T \mu E^2, 
        \end{aligned}
    \end{equation}
    for any $E \leq r_q$, any $x \in X_{\bar{V}_q}$ and any $T \leq T_1$.
    \qed
    \end{pf}
\end{prop}
\begin{defn}
    Define the following set:
    \begin{equation}
        \Sigma \vcentcolon = \{(x, z) \vcentcolon V(x) \leq \bar{V}_q, \| z - \bar{z}(x)\| \leq \tilde{r}_q\}.
    \end{equation}
\end{defn}
The following theorem shows positive invariance of $\Sigma$.
\begin{lem}[Invariance of $\Sigma$]\label{lem:invariance_q2}
    Let Assumptions \ref{assum:cont}, 
    \ref{assum:lyapunov_stability}, \ref{assum:sog},  
    \ref{assum:sreg}, \ref{assum:contraction_quad} and \ref{assum:state_evolution_bound} 
    hold. 
    Define 
    \begin{equation}\label{eq:ts_tilde_r_q}
        \tilde{r}_q \vcentcolon = \min \bigg\{ r_q, \sqrt{\frac{\bar{a}\bar{V}_q}{\mu}}\bigg\} \, \text{and} \, T'_2 \vcentcolon =\frac{(1-\kappa)\tilde{r}_q \sqrt{a_1}}{\sqrt{\bar{V}}\gamma}.
    \end{equation}
    Then, for any $(x,z) \in \Sigma$ and any $T \leq T_2 \vcentcolon = \min\{ T'_2, T_1\}$, it holds that 
    $(x_+, z_+) \in \Sigma$. Moreover, the following 
    coupled system-optimizer contractions hold:
    \begin{equation}\label{eq:coupled_contr}
    \begin{aligned}
        && V_+ & \leq (1 - T \bar{a}) V + T \mu E^2, \\
    && E_+    & \leq T \hat{\gamma} V^{\frac{1}{2}}  + \kappa  E, \\
    \end{aligned}
    \end{equation}
    where $\hat{\gamma} = \frac{\gamma}{\sqrt{a_1}}$.
    \begin{pf}
    Given that $E \leq \tilde{r}_q \leq r_q$ and $x \in X_{\bar{V}_q}$, 
    we can apply the contraction 
    from Proposition \ref{prop:inexact_lyapunov_contraction_q2}, such that 
    \begin{equation}
        V_+  \leq (1- T \bar{a}) V + T \mu E^2, 
    \end{equation}
    holds. Moreover, due to the definition of $\tilde{r}_q$, 
    $V_+ \leq \bar{V}_q$ holds, which implies that $x_+$ is in $X_{\bar{V}_q}$. Similarly, 
    due to the fact that $E \leq \tilde{r}_q \leq r_z$ and $x \in X_{\bar{V}_q} \subseteq X_{\bar{V}}$,
    we can apply the result from Proposition \ref{prop:contraction_x}, which shows that 
    \begin{equation}\label{eq:lem_contr_eq_1}
    E_+  \leq \, \kappa E + T \gamma \| x \|\quad \text{and} \quad E_+  \leq r_z
    \end{equation}
    must hold. 
    Using Assumption \ref{assum:lyapunov_stability} in Equation \eqref{eq:lem_contr_eq_1},
    we obtain
    \begin{equation}
        E_+  \leq \, \kappa E + T \hat{\gamma} V^{\frac{1}{2}}.
    \end{equation}
    Moreover, due to \eqref{eq:ts_tilde_r_q}, $E_+ \leq \tilde{r}_q$ holds.
    \qed
    \end{pf}
\end{lem}
Lemma \ref{lem:invariance_q2} shows that, we can guarantee that the state of the combined system-optimizer dynamics $(x,z)$ will 
not leave $\Sigma$ under the assumption that the sampling time $T$ is short enough.
Moreover, due to subadditivity of the square root, the following holds:
\begin{equation}
    V_+^{\frac{1}{2}} \leq (1 - T \bar{a})^{\frac{1}{2}} V^{\frac{1}{2}} + (T \mu)^{\frac{1}{2}}E
\end{equation}
such that we can regard the following simpler dynamics:
\begin{defn}[Auxiliary Dynamics]\label{defn:aux_dyn_q_2}
    We will refer to the following (linear) dynamical system:
\begin{equation}\label{eq:aux_dyn_q_2}
    \begin{aligned}
        &&\nu_+ &= (1 - T \bar{a})^{\frac{1}{2}} \nu + (T \mu)^{\frac{1}{2}} \epsilon, \\
        && \epsilon_+    & = T \hat{\gamma} \nu  + \kappa  \epsilon \\
    \end{aligned}
\end{equation}
with states $\nu, \epsilon \in \mathbb{R}$ as \textit{auxiliary dynamics}.
\end{defn}
\begin{rem}
    Notice that the considerations made by \cite{Diehl2007b}, in a similar setting, lead to the same type of coupled contraction from Lemma \ref{lem:invariance_q2}. 
    An attractivity proof that implicitly uses auxiliary dynamics that would be obtained directly from \eqref{eq:coupled_contr}
     is derived in \citep{Diehl2007b}. However, due to the fact 
    that such auxiliary system-optimizer dynamics are not Lipschitz at $(0,0)$, 
    it would not be possible to prove stability with standard linear analysis tools. 
\end{rem}
\section{Asymptotic Stability Result}
Due to linearity of the auxiliary dynamics \eqref{eq:aux_dyn_q_2}, 
we can study asymptotic stability with standard tools from linear systems analysis. 
\begin{thm}[Asymptotic Stability] \label{thm:as_stability}
    Let Assumptions \ref{assum:cont}, 
    \ref{assum:lyapunov_stability}, \ref{assum:sog}, \ref{assum:sreg}, \ref{assum:contraction_quad} 
    and \ref{assum:state_evolution_bound} hold. Then, the origin $(\nu, \epsilon) = (0,0)$ is asymptotically stable for the auxiliary dynamics \eqref{eq:aux_dyn_q_2}.
    \begin{pf}
        A sufficient and necessary condition for the asymptotic stability of \eqref{eq:aux_dyn_q_2} is that the eigenvalues 
of the matrix 
\begin{equation}
    A =
    \begin{bmatrix}
        (1 - T \bar{a})^{\frac{1}{2}} &\quad (T \mu)^{\frac{1}{2}} \\ 
        T \hat{\gamma}   &\quad \kappa
    \end{bmatrix},
\end{equation}
are smaller than one in absolute value. In order to compute the eigenvalues, we need to solve
\begin{equation*}
    \det (\lambda \mathbb{I} - A) = (\lambda - (1 - T \bar{a})^\frac{1}{2}) \cdot (\lambda - \kappa) + O\left(T^{\frac{3}{2}}\right) = 0, 
\end{equation*}
which entails
\begin{equation}
        \lambda_1 =(1-T\bar{a})^{\frac{1}{2}} + O\left(T^{\frac{3}{2}}\right),\, \lambda_2 =\kappa + O\left(T^{\frac{3}{2}}\right).
\end{equation}
Hence, for a sufficiently small sampling time $T$, the origin $(\nu, \epsilon) = (0,0)$ is asymptotically stable. \qed
\end{pf}
\end{thm}
\begin{figure}
    \centering
    \includegraphics[scale=0.7]{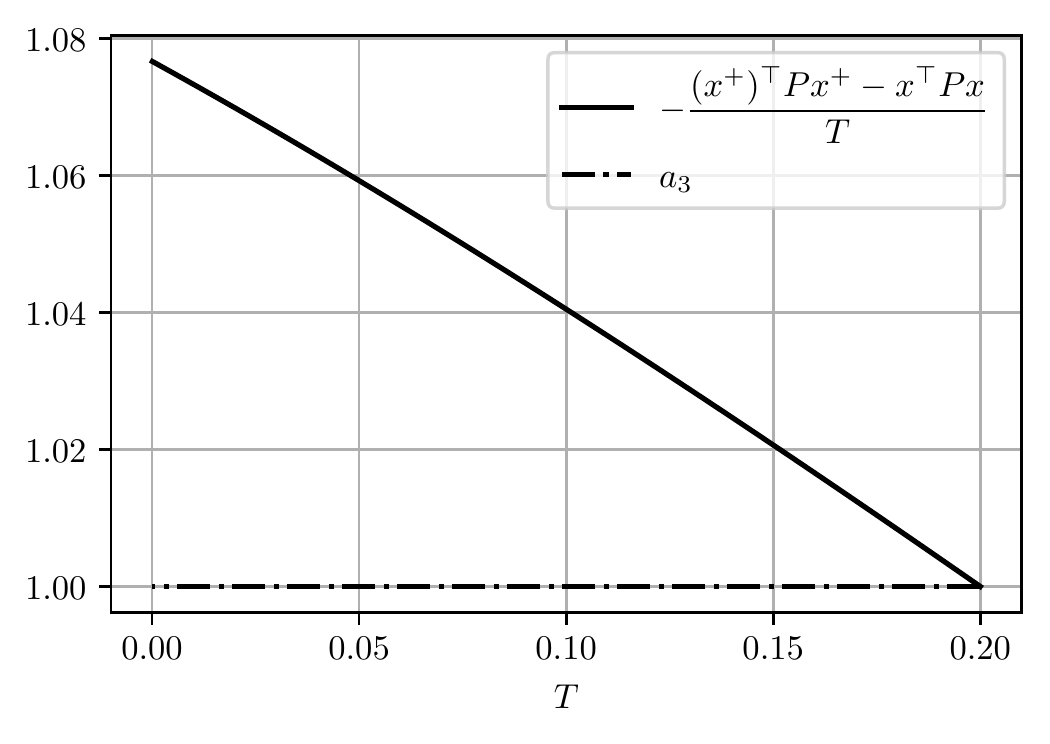}
    \caption{Lower bound on $V(x)$ and Lyapunov decrease as a function of $T$. Here we regard $x^+ \vcentcolon = (A_T + B_TK_d)x$ and check that $x^{\top}Px$ is 
    a still a Lyapunov function for such closed-loop system and estimate the decrease rate.}\label{fig:a_1}
\end{figure}
Theorem \ref{thm:as_stability} shows asymptotic stability of the auxiliary dynamics \eqref{eq:aux_dyn_q_2}, but not of the original system-optimizer dynamics \eqref{eq:so_dyn_compact}.
The following corollary shows that we can easily extend the stability result in this sense.
\begin{cor} Let Assumptions \ref{assum:cont}, 
    \ref{assum:lyapunov_stability}, \ref{assum:sog}, 
    \ref{assum:sreg}, \ref{assum:contraction_quad} and  \ref{assum:state_evolution_bound} hold. Then, the origin $\xi = (x, z) = (0,0)$ is locally asymptotically stable for the system-optimizer dynamics \eqref{eq:so_dyn_compact}.
    \begin{pf}
        Let $\chi \vcentcolon = (\nu, \epsilon)$ denote the state of the auxiliary dynamics in compact form and regard 
        the sequences $\{\chi_k\}$ and $\{\xi_k\}$ generated by the auxiliary and original coupled dynamics for any
        compatible initial conditions $\chi_0 = (\nu_0, \epsilon_0)$ and $\xi_0 = (x_0, z_0)$, with $V(x_0)^{\frac{1}{2}} = \nu_0$ and 
        $\|z_0 - \bar{z}(x_0)\| = \epsilon_0$, respectively. Due to the definition of the auxiliary dynamics, we have that, 
        for any compatible initial conditions chosen in such a way, $V(x_k)^{\frac{1}{2}} \leq \nu_k$ and $\| z_k - \bar{z}(x_k)\|\leq \epsilon_k$ 
        holds for any $k\geq 0$. Then, for any $k \geq 0$, we can write
        \begin{equation*}
        \begin{aligned}
            &&\| \xi_k \| &= \sqrt{\| x_k\|^2 + \| z_k - \bar{z}(x_k) \|^2} \leq \sqrt{\frac{1}{a_1}V(x_k) + E_k^2} \\
            && &\leq \sqrt{\frac{1}{a_1} \nu_k^2 + \epsilon_k^2} \leq \tilde{a}_1^{-\frac{1}{2}} \sqrt{\nu_k^2 + \epsilon_k^2} = \tilde{a}_1^{-\frac{1}{2}} \| \chi_k \|,
        \end{aligned}
        \end{equation*}
        where $\tilde{a}_1 = \min \{ a_1, 1\}$.

        Since $(\nu, \epsilon) = (0,0)$ is locally asymptotically stable for the auxiliary dynamics, we have that, for any $\epsilon' > 0$, there 
        exists a $\delta' > 0$ such that, if $\norm{\chi_0} < \epsilon'$, then $\norm{\chi_k} < \delta'$, for any $k \geq 0$.
        Since 
        \begin{equation}
            \| \chi_0 \|  \leq \sqrt{\tilde{a}_2} \norm{\xi_0},
        \end{equation}
        where $\tilde{a}_2 = \max \{ a_2, 1\}$, we can choose an arbitrarily small $\epsilon = \tilde{a}_2^{-\frac{1}{2}}\epsilon'$, 
        such that $\| \xi_0 \| \leq \epsilon = \tilde{a}_2^{-\frac{1}{2}}\epsilon' \implies \| \chi_0 \| \leq \epsilon' \implies \| \chi_k \| \leq \delta', \forall k \geq 0$ and  $\| \xi_k \| \leq \delta \vcentcolon = \tilde{a}_1^{-\frac{1}{2}}\delta', \forall k \geq 0$, which proves stability.
        \par
        Finally, local attractivity can be trivially shown by observing that $\lim_{k\to\infty} \| \chi_k \| = 0 \implies \lim_{k\to\infty} \| \xi_k \| = 0$.
        \qed
    \end{pf}
\end{cor}
\section{Illustrative Example}
In this section, although the results derived apply to a much more general class of problems (twice-continuously nonlinear dynamics and cost),
we discuss an illustrative numerical example where we exploit a simplified 
setting in order to be able to explicitly compute all the constants used in the assumptions of Theorem \ref{thm:as_stability}. 
In particular, we regard the following unconstrained, linear-quadratic optimal control problem:
\begin{equation}\label{eq:lqmpc}
\begin{aligned}
&\underset{\begin{subarray}{c}
    s(\cdot), u(\cdot)
\end{subarray}}{\min}	&&\int_{0}^{\infty} \begin{bmatrix}s(t) \\ u(t)\end{bmatrix}^{\top}\begin{bmatrix}Q_c & 0\\ 0& R_c\end{bmatrix}\begin{bmatrix}s(t) \\ u(t)\end{bmatrix} \\
                            &\quad \text{s.t.} 		   	&&s(0) - x = 0,\\
                            & 							&&\dot{s}(t) = A_c s(t) + B_c u(t), \quad  t \in [0,\infty], \\
\end{aligned}
\end{equation}
where the continuous-time dynamics are defined by
\begin{equation}
    A_c \vcentcolon= \begin{bmatrix}0 & \,\,\,1 \\ 0 &\,\,\, 0\end{bmatrix}, \quad B_c \vcentcolon= \begin{bmatrix}0 \\ 1 \end{bmatrix}
\end{equation}
and the matrices   
\begin{equation}
    Q_c \vcentcolon= \mathbb{I}_2 \quad \text{and}  \quad R_c \vcentcolon= 1
\end{equation}
 define the cost. 
\par
Problem \eqref{eq:lqmpc} is discretized using multiple shooting 
with a single shooting node and a fixed discretization time $T_d = 0.1 s$ as follows:
\begin{equation}\label{eq:lqmpc}
\begin{aligned}
&\underset{\begin{subarray}{c}
    s_0, s_1, u_0
\end{subarray}}{\min}	&&T_d \begin{bmatrix}s_0 \\ u_0\end{bmatrix}^{\top}\begin{bmatrix}Q_c & 0\\ 0& R_c\end{bmatrix}\begin{bmatrix}s_0 \\ u_0\end{bmatrix} + s_1^{\top}P s_1\\
                            &\quad \text{s.t.} 		   	&&s_0 - x = 0,\\
                            & 							&&s_1 = A_{T_d} s_0 + B_{T_d} u_0,
\end{aligned}
\end{equation}
where the discrete-time dynamics are obtained using an exact discretization with piece-wise constant 
parametrization of the control trajectories:
\begin{equation}
    A_{T_d} \vcentcolon= \exp{(A_c T_d)}, \, B_{T_d} \vcentcolon= \left(\int_{0}^{T_d}\exp{(A_c\tau) \text{d}\tau}\right)B_c.
\end{equation}
Finally, the symmetric positive-definite 
matrix $P$ that defines the terminal cost for the discretized 
problem is computed by solving the discrete-time algebraic Riccati equation
\begin{figure}
    \centering
    \includegraphics[scale=0.75]{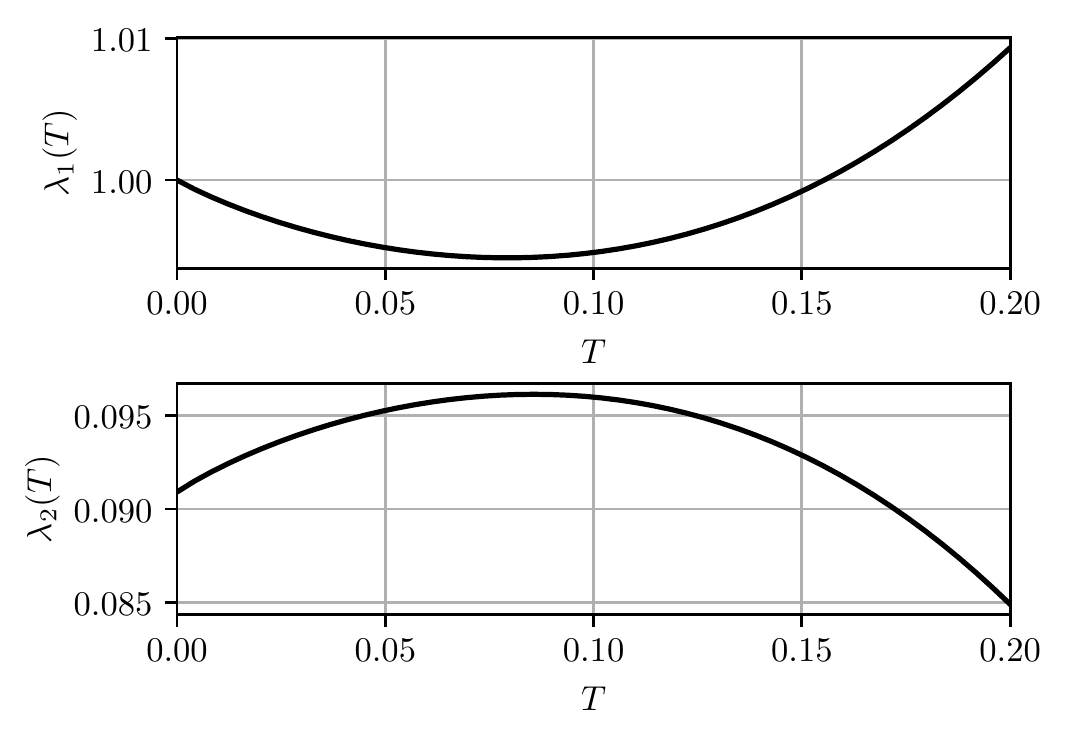}
    \caption{Eigenvalues as a function of the sampling time $T$ for problem \eqref{eq:lqmpc}. For sufficiently 
    short sampling times, the auxiliary system in Definition \ref{defn:aux_dyn_q_2} is asymptotically stable.}\label{fig:eigvals}
\end{figure}
\begin{equation}
    \begin{aligned}
        &P &&= A_{T_d}^{\top}PA_{T_d} - \\
        &  &&(A_{T_d}^{\top}PB_{T_d})(R + B_{T_d}^{\top}PB_{T_d})^{-1}(B_{T_d}^{\top}PA_{T_d}) + Q,
    \end{aligned}
\end{equation}
where $R \vcentcolon= T_d R_c$ and $Q \vcentcolon= T_d Q_c$.
After elimination of $s_0$, the first-order optimality conditions of problem \eqref{eq:lqmpc} read
\begin{equation}\label{eq:lqmpc_foc}
    H u_0 + Gx = 0,
\end{equation}
where 
\begin{equation}
    H \vcentcolon = (T_dR + B_{T_d}^{\top}PB_{T_d}), \quad G \vcentcolon=B_{T_d}^{\top}PA_{T_d}.
\end{equation}
We solve \eqref{eq:lqmpc_foc} with the following real-time gradient descent method:
\begin{equation}
    u_{0,+} = -\tilde{H}^{-1}\left((H-\tilde{H})u_0 + Gx_+\right),
\end{equation}
where $\tilde{H} = \rho \,\mathbb{I}$ for some positive constant $\rho > 1$. 
Using standard arguments from convergence theory for Newton-type methods (see e.g \citep{Diehl2016}), 
it is easy to show, that, for a fixed value of the parameter $x$, the following contraction estimate holds:
\begin{equation}
    \| u_{0,+} - \bar{u}_0(x) \| \leq \hat{\kappa} \| u_0 - \bar{u}_0(x) \|,  
\end{equation}
where 
\begin{equation}
    \hat{\kappa}\vcentcolon = \| \tilde{H}^{-1} (H-\tilde{H}) \|.
\end{equation}
Since $V(x) = x^{\top}Px$ and $\bar{u}_0(x) = -H^{-1}Gx$, we 
can compute exactly the constants $\mu = 2\lambda_{\max}(H)$ and $\sigma = \| H^{-1}G \|$. 
Let 
\begin{equation}
    A_{T} \vcentcolon= \exp{(A_c T)}, \, B_{T} \vcentcolon= \left(\int_{\tau = 0}^{T}\exp{(A_c\tau) \text{d}\tau}\right)B_c
\end{equation}
describe the dynamics discretized according to the sampling time $T$. Given that 
\begin{equation}
    x_+ = A_{T}x + B_{T}u= x + \frac{T}{T}((A_{T} - \mathbb{I})x + B_{T}u),
\end{equation}
we can compute the constants 
\begin{equation}
    L_{\psi, x} = \frac{1}{T}\| A_{T} - \mathbb{I}\| \quad \text{and} \quad L_{\psi, u} = \frac{1}{T} \| B_{T} \|.
\end{equation}
Following the definitions in Propositions \ref{prop:state_evolution_bound} and \ref{prop:contraction_x},
we have $\hat{\gamma} = \frac{1}{\sqrt{a_1}}\sigma\hat{\kappa}\eta$ and $\kappa = \hat{\kappa}(1+T \sigma \theta)$,
where $\eta = L_{\psi,u} + L_{\psi,x} \sigma$ and $\theta = L_{\psi, u}$. 
In order to validate Assumption \ref{assum:lyapunov_stability} and compute an estimate for constant $a_3$, we compute 
the largest eigenvalue $\lambda_{\max} (\Delta P)$ of the matrix 
\begin{equation}
    \Delta P(T) \vcentcolon = \frac{1}{T} ((x^{+})^{\top}Px^{+} - x^{\top}Px),
\end{equation}
where 
\begin{equation}
    x^{+} = (A_{T} + B_{T}K)x
\end{equation}
and $K = -H^{-1}G$. Figure \ref{fig:a_1} shows the estimated decrease rate compared with the minimum eigenvalue 
of $P$. Choosing 
\begin{equation}
    a_3 = \underset{T}{\min} \,\lambda_{\max} (\Delta P (T)),
\end{equation}
we obtain a value of $a_3$ that satisfies \eqref{eq:assum_lyap_2} for any $T$ such that $0 \leq T \leq T_d$.
Finally, we can compute the constants in \eqref{eq:assum_lyap_1} as $a_1 = \lambda_{\min} (P)$ and $a_2 = \lambda_{\max} (P)$.

Given that we can numerically compute all the constants involved in the Assumptions of Theorem \ref{thm:as_stability}, it 
is possible to compute the longest sampling time for which the auxiliary system in \eqref{eq:aux_dyn_q_2} is asymptotically stable. 
Figure \ref{fig:eigvals} shows the eigenvalues of the auxiliary system as a function of $T$.  
\section{Conclusions and Outlook}
In this paper we present an asymptotic stability proof for real-time methods for unconstrained NMPC.
We extend the well-known attractivity results derived in \citep{Diehl2005b} and \citep{Diehl2007b} under similar 
settings. The interaction between system and optimizer is analyzed using an auxiliary system whose states are the square root of the optimal value and the numerical error associated with the approximate solution. Asymptotic stability of such 
auxiliary system is established under the assumption that the sampling time is sufficiently short. 
Ongoing research involves the extension of the results to the case where the optimal control formulation includes inequality constraints 
and the derivation of practical conditions under which Assumptions \ref{assum:lyapunov_stability} and \ref{assum:sog} hold, using arguments 
along the lines of Remarks \ref{rem:d2c_lyap} and \ref{rem:sog}. Similarly, Assumption 
\ref{assum:sreg} will be adapted to the more general setting described in Remark 
\ref{rem:sreg}.
\bibliography{syscop} 
\end{document}